\documentclass{amsart}

\usepackage{amsmath}
\usepackage{amssymb}
\usepackage{amsfonts}
    
\newtheorem{theorem}{Theorem}

\newtheorem{lemma}[theorem]{Lemma}
\newtheorem{corollary}[theorem]{Corollary}
\newtheorem{remark}[theorem]{Remark}
\newtheorem{example}{Example}
\theoremstyle{definition}
\newtheorem{definition}[theorem]{Definition}
\theoremstyle{remark}

\begin{document}

\title{Genus  2 curves with $(3,3)$-split Jacobian and \\
large automorphism group}

\author{T. Shaska}

\address{University of California at  Irvine,  
Irvine, CA  92697.
}

\email{tshaska@math.uci.edu}

%***************************************
% my definitions
%****************************************
\def\bC{{\mathbb C}}
\def\C{{\mathcal C}}
\def\bQ{{\mathbb Q}}
\def\bZ{{\mathbb Z}}
\def\bP{{\mathbb P}}
\def\L{{\mathcal L}}
\def\M{{\mathcal M}}
\def\cO{{\mathcal O}}
\def\l{{\lambda}}
\def\sem{{\rtimes}}
\def\iso{{\, \cong\, }}
\def\p{\mathfrak p}
\def\E{\mathcal E}
\def\u{\mathfrak u}
\def\v{\mathfrak v}
\def\w{\mathfrak w}
\def\th{\theta}
\def\t_h{\theta_1}
\def\X{\mathfrak X}
\def\Y{\mathfrak Y}

%***************************************************************
\begin{abstract}
Let $\C$ be a genus 2 curve defined over $k$, $char (k) =0$. 
If  $\C$ has a $(3,3)$-split Jacobian  then  
we show that the automorphism group  $Aut(\C)$  is isomorphic to one of 
the following:  $\bZ_2, V_4, D_8$, or $D_{12}$.  
There are exactly six $\bC$-isomorphism classes of 
genus two curves $\C$ with $Aut(\C)$ isomorphic to $D_8$ (resp., $D_{12}$). 
%We compute their absolute invariants $i_1, i_2, i_3$. 
We  show that exactly four (resp., three) of these classes  
with group $D_8$ (resp., $D_{12}$) have  representatives defined over $\bQ$.
We discuss some of these curves in detail and find their rational points.
\end{abstract}

\maketitle

%***********************************************************************
\section{Introduction}
%&&&&&&&&&&&&&&&&&&&&&&&&&&&&&&&&&&&&&&&&&&&&&&&&&&&&&&&&&&&&&&&

Let $\C$ be a genus 2 curve defined over an algebraically closed field $k$, 
of characteristic zero. We denote by $K:=k (\C)$  its function field and by
 $Aut(\C):=Aut(K/k)$ the automorphism group  of $\C$. Let $\psi: \C \to \E$ be
a degree $n$ maximal covering (i.e. does not factor through an isogeny)  to an
elliptic curve $\E$ defined over $k$. 
We say that $\C$ has a {\it degree n  elliptic subcover}.  
Degree $n$ elliptic subcovers occur in pairs. 
Let $(\E, \E')$ be such a pair. It is well known that there is an
isogeny of degree $n^2$ between  the Jacobian $J_\C$ of $\C$
and  the product $\E \times \E'$. 
We say that $\C$ has {\bf (n,n)-split Jacobian}.
The locus of such $\C$ (denoted by $\L_n$) is  an algebraic 
subvariety  of the moduli space $\M_2$ of genus two curves. 
For the  equation of $\L_2$ in
terms of Igusa invariants, see \cite{Sh-V}. Computation of the equation of
 $\L_3$  was the main focus of \cite{Sh2}. For $n > 3$, equations
of $\L_n$ have not yet been computed. 

\par Equivalence classes of degree 2 coverings $\psi: \C \to \E$ are in  1-1 
correspondence with  conjugacy classes of non-hyperelliptic involutions in $Aut(\C)$. 
In any characteristic different from 2,  the automorphism group
  $Aut(\C)$  is isomorphic to one of the following: 
\ $\bZ_2$, $\bZ_{10}$, $V_4$, $D_8$, $D_{12}$, $\bZ_3 \sem D_8$, $ GL_2(3)$, 
or $2^+S_5$; see \cite{Sh-V}.
Here  $V_4$ is the Klein 4-group, $D_8$ 
(resp., $D_{12}$) denotes the dihedral group of order 8 (resp., 12),  and
$\bZ_2, \bZ_{10}$ are cyclic groups of order 2 and 10.
For  a description of other groups, see \cite{Sh-V}.
 If $Aut(\C)\iso \bZ_{10}$ then $\C$ is isomorphic 
to $Y^2=X^6-X$. Thus, if $\C $ has extra automorphisms and it is not isomorphic 
to $Y^2=X^6-X$ then $\C\in \L_2$. We say that  a genus 2 curve $\C$ has {\bf large
automorphism group} if the order of $Aut (\C)$ is bigger then 4. 
 
\par In section 2,  we describe the loci for genus 2 curves with $Aut(\C)$
isomorphic to $ D_8$ or $D_{12}$ in terms of Igusa invariants. 
From these invariants 
we are able to determine the field of definition of a curve $\C$ with 
$Aut(\C)\iso D_8$ or $D_{12}$. 
Further, we find the equation for this $\C$ and 
$j$-invariants of degree 2  elliptic subcovers in terms of $\, \, \, i_1, i_2, i_3$ 
(cf. section 2). 
This determines the fields of definition for these elliptic subcovers. 

\par Let $\C$ be a genus 2 curve with $(3,3)$-split Jacobian. In 
section 3 we give a short description of the space $\L_3$. 
Results described in section 3 follow from \cite{Sh2}, even though sometimes
nontrivially.  
We find equations of degree 3 elliptic subcovers
in terms of the coefficients of $\C$. 
 In section 4,  we show that  $Aut (\C)$  is one of the following: 
$\bZ_2$,  $V_4$, $D_8$, or  $D_{12}$.  Moreover,
 we show  that there are exactly six $\bC$-isomorphism classes of
genus two curves $\, \C\in \L_3$ 
with automorphism group $D_8$ (resp., $D_{12}$). 
We explicitly find the absolute invariants $i_1, i_2, i_3$  which determine 
these classes. 
For each such class we give the equation of a representative genus 2 curve $\C$. 
We notice that there are four cases (resp., three) such that the triple of 
invariants $(i_1,i_2, i_3)\in \bQ^3$ when $Aut(C)\iso D_8$ 
(resp., $Aut(C)\iso D_8$ ). 
Using results from section 2,  we determine that there are exactly 
four (resp., three)   genus 2  curves $\C\in \L_3$ 
(up to $\bar \bQ$-isomorphism)  with group $D_8$ (resp., $D_{12}$)
 defined over $\bQ$ and list their  equations  in Table 1.
We discuss these curves and their degree 2 and 3 
elliptic subcovers  in more detail in section 5. 
Our focus is  on the cases which have elliptic subcovers defined over $\bQ$.
In some of these  cases we are able to use these subcovers to
find  the rational points  of the genus 2 curve. This technique has been used 
before by Flynn and Wetherell \cite{FW} for degree 2 elliptic subcovers.

\par Curves of genus 2 with degree 2 elliptic subcovers (or with elliptic involutions)
were first studied by Legendre and Jacobi. The genus 2 curve with the largest
known number of rational points has automorphism group isomorphic to $D_{12}$; thus
it has degree 2 elliptic subcovers. It was found by Keller and Kulesz and it is known 
to have at least 588 rational points; see \cite{KK}.
Using degree 2 elliptic subcovers Howe, Leprevost, and Poonen 
\cite{HLP} were able to construct a family of  genus 2 curves
whose Jacobians each have large  rational torsion subgroups. 
Similar techniques probably could be applied using degree 3 elliptic subcovers. 
Curves of genus 2 with degree 3 elliptic subcovers
 have already occurred in the work of Clebsch,
 Hermite, Goursat, Burkhardt, Brioschi, and Bolza in the context of
elliptic integrals. For a history of this topic  see Krazer \cite{Krazer} (p.
 479). For more recent work see Kuhn \cite{Ku} and \cite{Sh2}.
More generally, degree $n$ elliptic subfields of genus 2 fields
have been studied by Frey \cite{Fr}, Frey and Kani \cite{FK}, Kuhn \cite{Ku},
and this author \cite{Sh1}. 

\par {\bf Acknowledgements:} The author wants to thank G. Cardona for pointing out 
the parametrizations in equations (5) and (6) which made some of 
the computations easier.

%\newpage
%&&&&&&&&&&&&&&&&&&&&&&&&&&&&&&&&&&&&&&&&&&&&&&&&&&&&&&&&&&&&&&&&&&&&&&&&&&&&
\section{Genus two curves  and  the moduli space $\M_2$.}
%&&&&&&&&&&&&&&&&&&&&&&&&&&&&&&&&&&&&&&&&&&&&&&&&&&&&&&&&&&&&&&&&&&&&&&&&&&&&
Let $k$ be an algebraically closed field of  characteristic  zero and 
$\C$ a genus 2 curve defined over $k$. Then $\C$ can be described as a double cover
of $\bP^1(k)$ ramified in 6 places $w_1, \dots , w_6$.
This sets up a bijection between isomorphism classes of genus 2 curves and
unordered distinct 6-tuples $w_1, \dots , w_6 \in \bP^1 (k)$ modulo automorphisms of
$\bP^1 (k) $.
An unordered 6-tuple $\{w_i\}_{i=1}^6$ can be described by a binary sextic,  
(i.e. a homogenous equation $f(X,Z)$ of degree 6). 
Let  $\M_2$ denote the moduli space of genus 2 curves; see \cite{Mu}. 
To describe $\M_2$ we
need to find polynomial functions of the coefficients of a binary sextic $f(X,Z)$ 
invariant under linear substitutions in $X,Z$ of determinant one. These invariants 
were worked out by Clebsch and  Bolza in the case of zero characteristic  and 
generalized by Igusa for any characteristic different from 2;
see \cite{Bo}, \cite{Ig}. 

%\subsection{Classical  invariants and the moduli space $\M_2$}
%***********************************************************************
\par Consider a binary sextic i.e. homogeneous polynomial
$f(X,Z)$ in $k[X,Z]$ of degree 6:
$$f(X,Z)=a_6 X^6+ a_5 X^5Z+\dots  +a_0 Z^6$$
{\it Classical  invariants} $\, \, \{ J_{2i} \}$ of $f(X,Z)$ are 
homogeneous polynomials of degree $2i$ in $k[a_0, \dots , a_6]$, 
for  $i=1,2,3,5$; see \cite{Ig}, \cite{Sh-V} for their 
definitions. 
Here  $J_{10}$ is simply the discriminant of $f(X,Z)$. 
It vanishes if and only if the binary sextic has a multiple linear factor. These 
$J_{2i}$    are invariant under the natural action of $SL_2(k)$ on
sextics. Dividing such an invariant by another one of the same degree
gives an invariant under $GL_2(k)$ action. 

\par Two genus  2 fields $K$ (resp., curves) in the standard form $Y^2=f(X,1)$ are
isomorphic if and only if the corresponding sextics are $GL_2(k)$
conjugate. Thus if $I$ is a $GL_2(k)$ invariant (resp., homogeneous $SL_2(k)$
invariant), then the expression $I(K)$ (resp., the condition $I(K)=0$)
is well defined. Thus the $GL_2(k)$ invariants are functions 
on the  moduli space $\mathcal M_2$ of genus 2 curves. This $\mathcal
M_2$ is an affine variety with coordinate ring 
$$k[\mathcal M_2]=k[a_0, \dots , a_6, J_{10}^{-1}]^{GL_2(k)}
={\textit{subring of degree 0 elements in}}$$
$k[J_2, \dots ,J_{10}, J_{10}^{-1}]$; see Igusa \cite{Ig}.
The  {\it absolute invariants}
\begin{equation}
i_1:=144 \frac {J_4} {J_2^2}, \quad i_2:=- 1728 \frac {J_2J_4-3J_6} {J_2^3}, \quad
i_3 :=486 \frac {J_{10}} {J_2^5}
\end{equation}
\noindent are even $GL_2(k)$-invariants.
Two genus 2 curves with $J_2\neq 0$ are isomorphic if and only if they have the same
absolute invariants. If  $J_2=0 $ then we can define new invariants as in  \cite{Sh2}. 
For the rest of this paper if we say 
``there is a genus 2 curve $\C$ defined over $k$'' 
we will mean the $k$-isomorphism class of $\C$.

\par One can define $GL_2(k)$  invariants with $J_{10}$ in the denominator which will 
be defined everywhere. However, this is not efficient in doing computations since 
the degrees of these rational functions in terms of the coefficients of $\C$ will be 
multiples of 10 and therefore  higher then degrees of $i_1, i_2, i_3$.
For the purposes of this paper defining $i_1, i_2, i_3$ as above 
 is not a restriction as it will be seen in 
the proof of the theorem 1.
For the proofs of the following two lemmas, see  \cite{Sh-V}.
%&&&&&&&&&&&&&&&&&&&&&&&&&&&&&&&&&&&

\begin{lemma}\label{thm1} The  automorphism group $G$ of a
genus 2 curve $\C$ in characteristic $\ne2$ is isomorphic to 
\ $\bZ_2$, $\bZ_{10}$, $V_4$, $D_8$, $D_{12}$, $\bZ_3 \sem D_8$, $ GL_2(3)$, 
or $2^+S_5$. 
The case when $G \iso 2^+S_5$ occurs only in characteristic 5. 
If $G \iso \bZ_3 \sem D_8$ (resp., $ GL_2(3)$)
then $\C$ has equation $Y^2=X^6-1$ (resp., $Y^2=X(X^4-1)$). If $G \iso \bZ_{10}$ then
$\C$ has equation $Y^2=X^6-X$.
\end{lemma}
\begin{remark} It is worth mentioning that the analogue of the above
lemma has been settled  for all $2 \leq g \leq 48$ in characteristic zero;  see \cite{MS}. 
This was a major computational effort using the computer  algebra system GAP. 
\end{remark}
For the rest of this paper we assume that $char(k)=0$. 
\begin{lemma}
i) The locus $\L_2$ of genus 2 curves $\C$ which have a degree 2 elliptic subcover
 is a closed subvariety of $\mathcal M_2$. The equation of $\L_2$ is
given by equation (17) in  \cite{Sh-V}.

ii) The locus  of genus 2 curves $\C$ with 
$Aut(\C)\iso D_8$ is given by the equation of $\L_2$  and 
\begin{small}
\begin{equation}
\label{D_8_locus}
\begin{split}
1706J_4^2J_2^2+2560J_4^3+27J_4J_2^4-81J_2^3J_6-14880J_2J_4J_
6+28800J_6^2 &=0
\end{split}
\end{equation}
\end{small}
ii) The locus  of genus 2 curves $\C$ with $Aut(\C)\iso D_{12}$ is 
\begin{small}
\begin{equation}
\label{D_12_locus}
\begin{split}
-J_4J_2^4+12J_2^3J_6-52J_4^2J_2^2+80J_4^3+960J_2J_4J_6-3600
J_6^2 &=0\\
864J_{10}J_2^5+3456000J_{10}J_4^2J_2-43200J_{10}J_4J_2^3-2332800000J_{10}^2-J_4^2J_2^6\\
-768J_4^4J_2^2+48J_4^3J_2^4+4096J_4^5 &=0\\
\end{split}
\end{equation}
\end{small}
\end{lemma}
%&&&&&&&&&&&&&&&&&&&&&&&&&&&&&&&&&
We we will refer to the locus of genus 2 curves $\C$ with $Aut(\C)\iso D_8$
 (resp., $Aut(\C)\iso D_8$ )   as $D_8$-locus (resp., $D_{12}$-locus).

\par Each genus 2 curve $\C\in \L_2$ has a non-hyperelliptic involution $v_0 \in Aut(\C)$. 
There is another non-hyperelliptic  involution $v_0^\prime:= \, v_0 \, w$,  
where $w$ is the hyperelliptic involution. 
Thus,  degree 2 elliptic subcovers come in pairs. 
We denote the pair of degree 2 elliptic subcovers by $(E_0, E_0^\prime)$.
If $Aut(\C)\iso D_8$ then $E_0 \iso E_0^\prime$ or $E_0$ and $E_0^\prime$ are
2-isogenous.  
If $Aut(\C)\iso D_{12}$, then $E_0$ and $E_0^\prime$ are isogenous of degree 3. 
See \cite{Sh-V} for details. 
\begin{lemma} Let $\C$ be a genus 2 curve defined over $k$. 
Then,

i)  $Aut(\C) \iso D_8$ if and only if  $\C$ is isomorphic to 
\begin{equation}\label{D_8}
Y^2=X^5+X^3+tX
\end{equation}
for some $t \in k \setminus \{ 0, \frac 1 4,  \, \frac {9} {100} \}$.

ii) $Aut(\C) \iso D_{12}$ if and only if  $\C$ is isomorphic to 
\begin{equation}\label{D_12}
Y^2=X^6+X^3+t
\end{equation}
for some $t \in k \setminus \{ 0, \frac 1 4, - \frac 1 {50} \}$.
\end{lemma}
\proof 
i) $Aut(\C)\iso D_8$: Then $\C$ is isomorphic to 
$$Y^2=(X^2-1)(X^4-\l X^2 +1)$$
for $\l \neq \pm 2$; see \cite{Sh-V}. 
Denote $\tau:=\sqrt{-2 \frac {\l+6} {l-2}}$.
The transformation 
$$\phi: (X,Y) \to (\frac {\tau x -1} {\tau x +1} ,  
\frac {4 \tau}  {(\tau x +1)^3} \cdot \frac {(\l+6)^2} {\l-2}       )$$ gives 
$$Y^2=X^5+X^3+tX$$
where $t= ( \frac {\l-2} { 2 (\l-6)} )^2  $ and $t\neq 0, \frac 1 4$. 
If $t=\frac 9 {100}$ then $Aut(\C) $ has order 24. 

\par Conversely, the absolute invariants $i_1, i_2, i_3$ of a genus  2 curve $\C$ 
isomorphic to 
$$Y^2=X^5+X^3+tX$$ 
satisfy the locus as described in lemma 2, part ii). Thus, $Aut(\C)\iso D_8$.

ii) $Aut(\C)\iso D_{12}$: In \cite{Sh-V} it  is shown that 
$\C$ is isomorphic to $$Y^2=(X^3-1)(X^3-\l)$$
for $\l\neq 0,1$ and $\l^2 -38\l +1 \neq 0$. Then,
$$\phi: (X,Y) \to  ( {(\l+1)}^{\frac 1 3} \, X, (\l +1) \, Y)  $$ 
transforms $\C$ to the curve with equation 
$$Y^2=X^6+X^3+ t $$
where $t=\frac \l {(\l+1)^2}$ and $t\neq 0, \frac 1 4$. If $t= - \frac 1 {50}$
then $Aut (\C)$ has order 48. 

\par The absolute invariants $i_1, i_2, i_3$ of a genus  2 curve $\C$ isomorphic to
$$Y^2=X^6+X^3+t$$ satisfy the locus as described in lemma 2, part iii). Thus,
$Aut(\C)\iso D_{12}$.  This completes the proof.

\qed

\par The following lemma determines a genus 2 curve for each point in the $D_8$ or 
$D_{12}$ locus.
\begin{lemma}  Let $\p:=(J_2, J_4, J_6, J_{10})$   be a  point in $\L_2$ such that 
$J_2\neq 0$  and $(i_1, i_2, i_3 )$ the corresponding absolute invariants. 

i) If $\p$ is in the $D_8$-locus, 
then the   genus two curve $\C$ corresponding to $\p$ is given by:
$$Y^2=X^5+X^3+
- \frac 3 4 \,  \frac {345i_1^2 + 50 i_1 i_2 - 90 i_2 -1296 i_1} 
{ 2925 i_1^2 + 250 i_1 i_2 -9450 i_2 -54000 i_1 + 139968} X$$
ii) If $\p$ is in the $D_8$-locus,
then the   genus two curve $\C$ corresponding to $\p$ is given by:
$$Y^2=X^6+X^4+
  \frac 1 4 \, \frac {540 i_1^2 + 100 i_1 i_2 -1728 i_1 +45 i_2} 
{2700 i_1^2 +1000 i_1 i_2 + 204525 i_1 + 40950 i_2 - 708588}$$
\end{lemma}
\proof
i) By previous lemma every 
genus 2 curve $\C$ with automorphism group $D_8$ is isomorphic to $Y^2=X^5+X^3+tX$. 
Since  $J_2\neq 0$ then $t\neq - \frac 3 {20}$ and  the  absolute invariants are:
\begin{small}
\begin{equation}
i_1 = \, -144\, t \, \frac {(20t-9)} {(20t+3)^2},  \quad
i_2 = \, 3456\, t^2\,  \frac {(140t-27)} {(20t+3)^3}, \quad
i_3 = \, 243\, t^3 \, \frac {(4t-1)^2} {(20t+3)^5}
\end{equation}
\end{small}
From the above system we have 
$$t= - \frac 3 4 \,  \frac {345i_1^2 + 50 i_1 i_2 - 90 i_2 -1296 i_1} 
{ 2925 i_1^2 + 250 i_1 i_2 -9450 i_2 -54000 i_1 + 139968}$$
ii) By previous lemma every 
genus 2 curve $\C$ with automorphism group $D_{12}$ is isomorphic to $Y^2=X^6+X^3+t$. 
The  absolute invariants are:
\begin{small}
\begin{equation}
i_1 = \, 1296 \frac {t(5t+1)} {(40t-1)^2},  \quad
i_2 = \, -11664 \frac {t(20t^2+26t-1)} {(40t-1)^3}, \quad
i_3 = \, \frac {729} {16} \, \frac {t^2 (4t-1)^3} {(40t-1)^5}  
\end{equation}
\end{small}
From the above system we have 
$$t=  \frac 1 4 \, \frac {540 i_1^2 + 100 i_1 i_2 -1728 i_1 +45 i_2} 
{2700 i_1^2 +1000 i_1 i_2 + 204525 i_1 + 40950 i_2 - 708588}$$
This completes the proof.

\qed

\noindent {\bf Note:} 
If $J_2=0$ then there is exactly one isomorphism class of genus 2 curves with
automorphism group $D_8$ (resp., $D_{12}$) given by $Y^2=X^5+X^3- \frac 3 {20} X$
(resp.,  $Y^2=X^6+X^3- \frac 1 {40} $).
\begin{remark} If the invariants $i_1, i_2, i_3\in \bQ$ then from above lemma there 
is a $\C$ corresponding to these invariants defined over $\bQ$. 
If a genus 2 curve does not have extra automorphisms (i.e. $Aut(\C)\iso \bZ_2$),  
then an algorithm of Mestre determines if the curve is defined over $\bQ$.
\end{remark}

If the order of the automorphism group 
$Aut (C) $ is divisible by 4, then $\C$ has degree 2 elliptic subcovers. 
These  elliptic subcovers are determined explicitly in \cite{Sh-V}. 
Do these elliptic subcovers of $\C$ have the same field of definition as $\C$? 
In general the answer is negative. 
The following lemma  determines the field of definition
of these elliptic subcovers when $Aut (\C)$ isomorphic to $ D_8$ or $D_{12}$.
%&&&&&&&&&&&&&&&&&&&&&&&&&&&&&&&&&&&&&&&&&&&&&&&
% lemma 5
%  ********************************** 
\begin{lemma} Let $\C$ be a genus 2 curve defined over $k$, $char(k)=0$.

i) If  $\C$ has  equation
$$Y^2=X^5+X^3+tX$$ 
where $t \in k\setminus \{ \neq 0, \frac 1 4, \frac 9 {100}\}$,
then its degree 2 elliptic subfields have j-invariants given by
$$j^2- 128 \frac {2000 t^2 + 1440 t +27} {(4t-1)^2 } + 
4096 \frac {(100t-9)^3} {(4t-1)^3}=0$$
ii) If  $\C$ has  equation
$$Y^2=X^6+X^3+t$$ 
where $t \in k\setminus \{ \neq 0, \frac 1 4, - \frac 1 {50}\}$,
then its degree 2 elliptic subfields have j-invariants given by
$$j^2- 13824\,  t \,  \frac {500 t^2 + 965 t +27} {(4t-1)^3 } + 
47775744\,  t \,  \frac {(25 t-4)^3} {(4t-1)^4}=0$$
\end{lemma}
\proof
The proof is elementary and follows from \cite{Sh-V}. \qed

%&&&&&&&&&&&&&&&&&&&&&&&&&&&&&&&&&&&&&&&&&&&&&&&&&&

%&&&&&&&&&&&&&&&&&&&&&&&&&&&&&&&&&&&&&&&&&&&&&&&&&&&&&&&&&&
\section{Curves of genus 2 with degree 3 elliptic subcovers}
%*********************************************************

%&&&&&&&&&&&&&&&&&&&&&&&&&&&&&&&
In this section we will give a brief description of spaces  $\L_2$ and $\L_3$. 
In case $J_2\neq 0$ we take these 
 spaces as  equations in terms of $i_1, i_2, i_3$, otherwise 
as  homogeneous equations in terms of $J_2, J_4, J_6, J_{10}$. By a point
$\p \in \L_3$ we will mean a tuple $(J_2, J_4, J_6, J_{10})$
which satisfies equation of $\L_3$. 
When it is clear that $J_2\neq 0$ then $\p\in \L_3$
 would mean  a  triple $(i_1,i_2,i_3)\in \L_3$. As before $k$ is an 
algebraically closed field of characteristic zero.
\begin{definition}
A {\bf non-degenerate pair} (resp., {\bf degenerate pair}) is a pair $(\C,\E)$ such
that $\C$ is a genus 2 field with a degree 3 elliptic subcover $\E$ where
$\psi: \C \to \E$ is ramified in two (resp., one) places.  Two such
pairs $(\C, \E)$ and $(\C' , \E')$ are called isomorphic if there is a
$k$-isomorphism $\C \to \C'$ mapping $\E\to \E'$.
\end{definition}
If $(\C, \E)$ is a non-degenerate pair, then $\C$ can be parameterized as follows 
\begin{equation}\label{eq_F1_F2}
Y^2=(\v^2 X^3+\u \v X^2 +\v X+1)\, (4\v^2 X^3 +\v^2 X^2+2\v X+1)
\end{equation}
where $\u, \v  \in k$ and the discriminant 
$$\Delta = -16\, \v^{17} \, (\v -27) \, 
(27 \v + 4 \v^2 - \u^2 \v + 4\u^3 -18 \u \v )^3 $$
of the sextic is nonzero. We let 
$R:=(27 \v + 4 \v^2 - \u^2 \v + 4\u^3 -18 \u \v )\neq 0$.
For $4\u - \v -9\neq 0$ degree 3 coverings are given by 
\begin{small}
$$U_1= \frac {\v X^2} {\v^2 X^3 + \u \v X^2 + \v X +1}, \quad 
U_2= \frac {(\v X +3)^2\, \,  ( (4\u-\v-9) X +3\u -v)} 
{\v \,(4\u -\v -9) (4 \v^2 X^3 + \v^2 X^2 + 2 \v X + 1) } 
$$
\end{small}
and elliptic curves have equations:
\begin{small}
\begin{equation}
\begin{split}
\E: & \quad V_1^2= R \, U_1^3 - (12 \u^2 - 2 \u \v - 18 \v )U_1^2 + 
(12\u -v) U_1 -4 \\
\E': & \quad V_2^2=c_3 U_2^3 +c_2 U_2^2 + c_1 U_2 +c_0
\end{split}
\end{equation}
\end{small}
where 
\begin{small}
\begin{equation}
\begin{split}
c_0 & = - (9\u -2 \v -27 )^3\\
c_1 & = (4\u -v -9) \, (729 \u^2 + 54 \u^2 \v -972\u \v - 18\u \v^2 +189\v^2 +
729 \v +\v^3) \\
c_2 & =- \v \,  (4\u -v -9)^2 \,  (54\u +\u \v -27\v)     \\
c_3 & =\v^2 \, (4\u -v -9)^3 \\
\end{split}
\end{equation}
\end{small}
The above facts can be deducted from lemma 1 of \cite{Sh2}. The case 
$4\u - \v -9= 0$ is treated separately in \cite{Sh2}.
There  is an   automorphism $\, \, \beta \in Gal_{k(u,\v)/k(i_1, i_2, i_3)}$ 
\begin{small}
\begin{equation}\label{eq_nu}
\begin{split}
\beta(\u) & =\frac {(\v-3\u)(324\u^2+15\u^2\v-378\u\v-4\u\v^2+243\v+72\v^2)}
{(\v-27)(4\u^3+27\v-18\u\v-\u^2\v+4\v^2)}\\
\beta(\v) & =- \frac {4(\v-3\u)^3}{4\u^3+27\v-18\u\v-\u^2\v+4\v^2} \\
\end{split}
\end{equation}
\end{small}
which permutes the $j$-invariants of $\E$ and $\E'$. 
 The map $$\th: (\u,\v) \to (i_1, i_2, i_3)$$
defined when $J_2\neq 0$ and $\Delta \neq 0$   has degree 2. 
Denote by $J_\th$ the 
Jacobian matrix of   $\th$.  Then $det (J_\th)=0$ 
consist of the (non-singular) curve $\X$ given by
\begin{equation}
\begin{split}\label{n_3_iso1}
\X: \quad 8\v^3+27\v^2-54\u \v^2-\u^2\v^2+108\u^2\v+4\u^3\v-108\u^3=0\\
\end{split}
\end{equation}
and  6 isolated $(\u, \v)$  solutions. These solutions correspond to the following
values for $(i_1, i_2, i_3)$:
%&&&&&&&&&&&&&&&&&&&&&&&&&&&&&&&&&&&&&&&&&&&&&&&&&&&&&&
\begin{scriptsize}
\begin{equation}
\label{exp_pts}
( - \frac {8019} {20}, -\frac {1240029} {200}, - \frac {531441} {100000}  ), \, \, 
 (\frac {729} {2116}, \frac {1240029}  {97336}, \frac {531441} {13181630464}), \, \,
(81, - \frac {5103} {25}, -\frac {729} {12500} )
\end{equation}
\end{scriptsize}
\noindent We denote  the image of $\X$ in the $\L_3$	locus by $\Y$. 
The map $\th$ restricted to $\X$ is unirational. The curve $\Y$ can be  computed as
an affine curve in terms of $i_1, i_2$. 
For each point $\p\in \Y$  degree 3 elliptic subcovers are isomorphic. 
If $\p$ is an ordinary point in $\Y$ and $\p \neq \p_6   $ (cf. Table 1) then 
the corresponding curve $\C_\p$ has automorphism group $V_4$.
%&&&&&&&&&&&&&&&&&&&&&&&&&&&&&&&&&&&
\par If $(\C, \E)$ is a degenerate pair then $\C$ can be parameterized as follows 
$$Y^2=(3X^2+4)(X^3+X+c)$$
for some $c$ such that $c^2\neq  - \frac 4 {27} $; see \cite{Sh2}.
We define $ \w :=c^2$. 
The map 
$$\w \to (i_1, i_2, i_3)$$
is injective  as it was shown in \cite{Sh2}. 
\begin{definition}
Let $\p$ be a point in $\L_3$. We say $\p$ is a {\bf generic 
point} in $\L_3$ if the corresponding 
$(\C_\p, \E)$ is a non-degenerate pair.
We define 
$$e_3 (\p):=
\left\{
\aligned
| \th^{-1} (\p)|,  \quad if \,\, \p \,\, is \, \, a \,\, generic \,\, point\\
1 \quad \quad \quad \quad otherwise \quad \quad \quad \quad  \quad  \, \, \,  \\
\endaligned
\right.
$$
\end{definition}
In \cite{Sh2} it is shown that 
the pairs $(\u,\v)$ with $\Delta (\u,\v)\neq 0$  bijectively parameterize the
isomorphism classes of non-degenerate pairs   $(\C,\E)$.
Those $\w$ with $\w\neq - \frac 4 {27}$ bijectively parameterize the
isomorphism classes of degenerate pairs   $(\C, \E)$.
Thus, the number $e_3(\p)$ is the number of 
isomorphism classes of such pairs   $(\C, \E)$. In \cite{Sh2} it is shown that 
$e_3(\p)=0,1,2$, or 4.
The following lemma describes the locus $\L_3$. For details see \cite{Sh2}.
\begin{lemma}
The locus $\L_3$ of genus 2 curves with degree 3 elliptic subcovers is the
closed subvariety of $\M_2$ defined by the equation 
\begin{equation}
\label{eq_L3_2}
C_8 J_{10}^8 + \dots  + C_1 J_{10} +C_0=0
\end{equation}
where coefficients $C_0, \dots ,C_8\in k[J_2, J_6, J_{10}]\, $ are displayed  in \cite{Sh2}. 
\end{lemma}
As noted above, with the assumption $J_2\neq 0$ equation \eqref{eq_L3_2}
can be written in terms of $i_1, i_2, i_3$. 
%***********************************************************
\section{Automorphism groups of genus  2 curves with degree 3 
elliptic subcovers}
%&&&&&&&&&&&&&&&&&&&&&&&&&&&&&&&&&&&&&&&&&&&&&&&&&&&&&&&&&
Let  $\C\in \L_3$ be a genus 2 curve defined over an algebraically closed
field  $k$, $char(k)=0$. 
The following theorem determines the automorphism group of $\C$. 
\begin{theorem}
Let $\C$ be a genus two curve which has a degree 3 elliptic subcover. 
Then the automorphism group
of $\C$ is one of the following: $\bZ_2, V_4$, $ D_8$, or $D_{12}$.
Moreover, there are exactly six curves $\C\in \L_3$ with automorphism
group $D_8$ and six curves $\C\in \L_3$ with automorphism group $D_{12}$. 
\end{theorem}
\proof
We denote by $G:=Aut(\C)$. None of the curves 
$Y^2=X^6-X$, $Y^2=X^6-1$, $Y^2=X^5-X$ have degree  3 elliptic subcovers since
their $J_2, J_4, J_6, J_{10}$ invariants don't satisfy equation \eqref{eq_L3_2}. 
From lemma 1 we have the following cases:

i) If $G\iso D_8$, then $\C$ is isomorphic to $$Y^2=X^5+X^3+t\,X$$
as in Lemma 3. 
Classical   invariants are:
\begin{small}
$$
J_2 = 40t+6, \, \, J_4  =4t(9-20t), \,\,  J_6  =8t(22t+9-40t^2), \, \,
J_{10} =16t^3(4t-1)^2.  
$$
\end{small}
\noindent Substituting in the equation \eqref{eq_L3_2} we have the following equation:
\begin{small} \begin{equation}
 (196t-81)^4 (49t-12) (5t-1)^4 (700t+81)^4 (490000\, t^2-136200 \, t+ 2401)^2 =0
\end{equation}  \end{small}
\noindent For 
$$\, t= \frac {81} {196}, \frac {12} {49}, \frac 1 5 , - \frac {81} {700}\, $$
the triple $(i_1, i_2, i_3)$ has the following values respectively: 
\begin{small}
$$
(\frac {729} {2116},  \frac {1240029}  {97336},\frac {531441} {13181630464}), \quad 
(\frac {4288} {1849}, \frac {243712} {79507}, \frac {64} {1323075987}),
$$
$$
(\frac {144} {49}, \frac {3456} {8575}, \frac {243} {52521875}), \quad
(-\frac {8019} {20}, -\frac {1240029}  {200},  -\frac {531441} {10000} )
$$
\end{small}
\noindent If $$490000\, t^2-136200 \, t+ 2401=0$$ then we have two distinct 
triples $(i_1, i_2, i_3 )$ which are in $ \bQ (\sqrt{2})$.
Thus, there are exactly 6 genus 2 curves 
$\C \in  \L_3$  with automorphism group $D_8$ and only four of them have rational
invariants.
%&&&&&&&&&&&&&& D_12 &&&&&&&&&&&&&&&&&&&&&&&&&&&&&&&&&&&&&&&

ii) If $G \iso D_{12}$ then $\C$ is isomorphic to a genus 2 curve in the
form $$Y^2=X^6+X^3+t$$
as in Lemma 3. 
Then, $J_2 =-6(40t-1)$ and 
\begin{small}
$$
J_4 =324t(5t+1), \, \, 
J_6 =-162t(740t^2+62t-1), \, \, 
J_{10} =-729t^2(4t-1) 
$$
\end{small}
Then the equation of $\L_3$ becomes:
\begin{small} 
\begin{equation}
(25t-4)\,(11t+4)^3\, (20t-1)^6\,(111320000t^3-60075600t^2+13037748t+15625)^3 =0
\end{equation} 
\end{small}
For $$t=\frac 4 {25}, - \frac 4 {11}, \frac 1 {20}$$  the corresponding values
for $(i_1,i_2,i_3)$ are respectively:
\begin{small} 
$$
(\frac {64} 5, \frac {1088} {25}, \frac 1 {84375}), \quad 
( \frac {576} {361}, \frac {60480} {6859}, \frac {243} {2476099}), \quad
(81, - \frac {5103} {25}, -\frac {729} {12500})
$$
\end{small}
If $$111320000t^3-60075600t^2=13037748t+15625=0$$ then there are three distinct 
triples $(i_1, i_2, i_3)$ none of which is rational.
Hence, there are exactly 6 classes of 
genus 2 curves $\C\in \L_3$ with $Aut(\C)\iso D_{12}$
of which three have rational invariants.

iii) $G\iso V_4$. There is a 1-dimensional family of genus 2 curves 
with a degree 3 elliptic subcover and automorphism group $V_4$ given by
$\Y$.

iv)  Generically genus 2  curves $\C $ have  $Aut(\C)\iso \bZ_2$. 
For example, every point  $\p \in \L_3 \setminus \L_2 $ correspond to a class 
of genus 2 curves  with degree 3 elliptic subcovers and  automorphism  
group  isomorphic to $\bZ_2$.
This completes the proof.

\qed

The theorem determines that are exactly 12 genus 2 curves $\C \in \L_3$ with
automorphism group $D_8$ or $D_{12}$. Only seven of them have rational invariants. 
From Lemma 4, we have the following:
\begin{corollary}
There are exactly four (resp., three) genus 2 curves $\C$ defined over 
$\bQ$ (up to $\bar \bQ$-isomorphism)
 with a degree 3 elliptic subcover which have automorphism group 
$D_8$ (resp., $D_{12})$. They are listed  in Table 1. 
\end{corollary}
\begin{scriptsize}
\begin{table}[!hbt]
\renewcommand\arraystretch{1.5}
\noindent\[
\begin{array}{|c|c|c|c|c|}
\hline
 &  \C & \p=(i_1, i_2, i_3) &e_3(\p) & Aut(\C) \\
\hline
\p_1 &196 X^5+ 196 X^3+ 81 X &
i_1=\frac {729} {2116},  i_2=\frac {1240029}  {97336},  i_3=\frac {531441} {13181630464}  &
2& D_8 \\
\hline
\p_2 &49 X^5+ 49 X^3+ 12 X &
i_1=\frac {4288} {1849}, i_2=\frac {243712} {79507},
i_3=\frac {64} {1323075987}  & 1& D_8 \\
\hline
\p_3 &5 X^5+ 5 X^3+  X & 
i_1=\frac {144} {49}, i_2=\frac {3456} {8575}, i_3=\frac {243} {52521875} & 2& D_8  \\
\hline
\p_4 & 700 X^5+ 700 X^3- 81  X &
i_1=-\frac {8019} {20}, i_2=- \frac {1240029}  {200}, i_3=- \frac {531441} {10000} &2& D_8    \\
\hline
\p_5 & 25 X^6+ 25 X^3 + 4 & 
i_1=\frac {64} 5, i_2=- \frac {1088} {25}, i_3=-\frac 1 {84375}   & 1& D_{12} \\
\hline
\p_6 & 11 X^6+ 11 X^3- 4 &
i_1= \frac {576} {361}, i_2=\frac {60480} {6859}, i_3=\frac
{243} {2476099}    & 1 & D_{12} \\
\hline
\p_7 & 20 X^6+20 X^3+ 1 & 
i_1=81, i_2=- \frac {5103} {25}, i_3=-\frac {729} {12500}   & 2& D_{12} \\
\hline
\end{array}
\]
\caption{Rational points  $\p \in \L_3$ with $|Aut(\p)| > 4$}
\label{table}
\end{table}
\end{scriptsize}
%
%\vspace{-.3in}
\begin{remark}
All points $\p$ in Table 1 are in the locus $det (J_\th ) =0$.
We have already seen cases  $\p_1, \p_4$, and  $\p_7$  as 
the exceptional points of $det (J_\th) =0$; see equation \eqref{exp_pts}. 
The class $\p_3$ is a  singular point of order 2  of $\Y$,
 $\p_2$ is the only point which  belong to  the degenerate case, and $\p_6$ is 
the only ordinary point in $\Y$ such that  the order of  
$Aut (\p)$ is greater then 4.
\end{remark}
%&&&&&&&&&&&&&&&&&&&&&&&&&&&&&&&&&&&&&&&&&&&&&&&&&&&&&&&&&&&
\section{Computing elliptic subcovers}
%&&&&&&&&&&&&&&&&&&&&&&&&&&&&&&&&&&&&&&&&&&&&&&&&&&&&&&&&&&
Next we will consider all points $\p$ in Table 1 and compute j-invariants 
of their degree 2 and 3 elliptic subcovers. To compute j-invariants of
degree 2 elliptic subcovers we use lemma 5 and the values of $t$ from the proof
of  theorem 1. We recall that for $\p_1, \dots , \p_4$ 
there are four degree 2 elliptic subcovers which are two and two isomorphic. 
We list the j-invariant of each isomorphic class. 
They are 2-isogenous as mentioned before. 
For $\p_5, \p_6, \p_7$ there are two degree 2 elliptic subcovers which are 
3-isogenous to each other. 
To compute  degree 3 elliptic subcovers for each
$\p$ we find the pairs $(\u,\v)$ in the fiber $\th^{-1}(\p)$ and 
then use equations (9).
We focus on cases which have elliptic  subcovers defined over $\bQ$. There are 
techniques  for  computing rational points of a genus two curves which have
degree 2 subcovers defined over $\bQ$ as in  Flynn and  Wetherell \cite{FW}. 
Sometimes degree 3 elliptic subcovers are defined over $\bQ$ even though 
degree 2 elliptic subcovers are not; see  examples 2 and 6. 
These degree 3 subcovers help determine rational points of genus 2 curves as
illustrated in examples 2, 4, 5, and 6. 

\begin{example}  $\p=\p_1$: 
The j-invariants   of degree  3 elliptic subcovers are  
$j=j^\prime=66^3$. 
A genus 2 curve $\C$ corresponding to $\p$ is 
$$\C: \, \, Y^2=X^6+3X^4-6X^2-8.$$
{\it Claim: The equation above has no rational solutions. }

\medskip
Indeed, two  of degree 2 elliptic subcovers (isomorphic to each other) have  equations
$$\E_1: \, \, Y^2=x^3+3x^2-6x-8$$
$$\E_2: \, \, Y^2=-8x^3 -6x^2+3x+1$$
\noindent where $x=X^2$ (i.e. $\phi: \C \to \E_1$ of degree 2 such that 
$\phi (X,Y)=(X^2,Y) \, $ ). 
The elliptic curve $\E_1$ has rank 0. Thus, the rational points of $\C$ are 
the preimages of the torsion points of $\E_1$. The torsion group of $\E_1$ 
has order 4 and is given by $$Tor (\E_1)=\{ \cO, (-1,0), (2,0), (-4,0) \}$$ 
None of the preimages is rational. 
Thus, $\C$ has no rational points except the point at infinity.

\end{example}
\begin{example} $\p=\p_2$: The 
j-invariants   of degree  2 elliptic subcovers are
$$\, 76771008\pm 44330496\sqrt{3}.$$
The point $\p_2$ belongs to the degenerate locus with $\w=0$.
Thus, the equation of the  genus 2 curve $\C$ corresponding to $\p$ is 
$$\C: \quad Y^2=(3\,X^2+4)\, (X^3+X).$$
Indeed, this curve has both pairs $(\C, \E)$ and $(\C, \E')$
degenerate pairs. 
It is the only such genus 2 curve defined over $\bQ$. This fact
was noted in \cite{Ku} and \cite{Sh1}. 
Both authors failed to identify the automorphism group. 
The degree 3 coverings are 
$$U_1=\phi_1(X)=X^3+X, \quad U_2=\phi_2(X)= \frac {X^3} {3X^2+4} $$
and equations of elliptic curves:
$$\E:  \, \, V_1^2=27 U_1^3+4U_1,\quad and  \quad \E': \,\, V_2^2=U_2^3+ U_2. $$
$\E$ and $\E'$ are isomorphic with  $j$-invariant 1728. They have rank 0 and rational
torsion group of order 2, $\, Tor(\E)=\{\cO, (0,0)\}$. Thus, the only rational points
of $\C$ are in fibers $\, \, \phi_1^{-1}(0)\, $ and $\, \, \phi_2^{-1}(\infty)$. 
Hence, $\C(\bQ)=\{(0,0), \infty \}$.
\end{example}
\begin{example} $\p=\p_3$:
All degree 2 and 3 elliptic subcovers are defined over $\bQ (\sqrt{5})$.
\end{example}
%
%*************************************************
\begin{example}  $\p=\p_4$: 
Degree 2 elliptic subcovers have j-invariants 
$$\frac {1728000} {2809} \pm \frac {17496000} {2809} \sqrt{I}  $$ where $I^2=-1$. 
Thus,
we can't recover any information  from degree 2 subcovers. 
One  corresponding value for  $(\u,\v)$ is  $( \frac {25} 2, \frac {250} 9)$. 
Then $\C$ is 
$$\C: \quad 3^8 \cdot Y^2= (100X+9)(2500X^2+400X+9) \, (25X+9)(2500X^2+225X+9)  $$
Degree 3 elliptic subcovers are have equations 
\begin{equation}
\begin{split}
\E :& \quad V_1^2= - \frac 1 {81} (10 U_1-3) (8575 U_1^2 -2940 U_1 +108) \\
\E' :& \quad  V_2^2= - \frac {686} {59049} \, (1700 U_2 -441)
(1445000 U_2^2-696150 U_2 + 83853)\\
\end{split}
\end{equation}
where 
\begin{small}
\begin{equation}
\begin{split}
U_1  & =\phi_1(X)=  2250 \frac {X^2} {(25X+9)(2500X^2+400X+9)} \\
U_2  & =\phi_2(X)=
  \frac 1 {340} \, \frac {(250X+27)^2 (340X+9)} {(100X+9)(2500X^2+400X+9)}\\
\end{split}
\end{equation}
\end{small}
The curve $\E'$ has rank 0 and torsion group 
$$Tor (\E') =\{ \cO, (\frac {441} {1700}, 0)\}.$$ 
Thus, the rational points of $\C$ are in $\phi_2^{-1}(\frac {441} {1700})$ 
and $\phi_2^{-1}(\infty)$.
Thus, $$\C(\bQ)=\{ (- \frac 9 {100}, 0 ), (- \frac 9 {25} ,0)\}.$$
\end{example}
%
%**************************************************
%
\begin{example} $\p=\p_5$: Degree 2  $\, j$-invariants are 
$j_1=0$, $ j_2=-1228800$ and degree 3 
$\, j$-invariants as shown below are $j=j'=0$.
Let $\C$ be  genus 2 curve with  equation
$$\C: \quad Y^2=(X^3+1)(4X^3+1)$$
corresponding to $\p$. 
The case  is treated separately in \cite{Sh2}.  Degree 3 elliptic subcovers
have equations 
$$\E: \, \,  V_1^2=-27 U_1^3 +4 , \quad \E': \, \,  V_2^2=-16 (27 U_2^3 -1)$$
where $$U_1= \phi_1 (X)= \frac {X^2} {X^3+1}, \quad U_2= \phi_2 (X) =\frac X {4X^3+1}.$$
The rank of both $\E$ and $\E'$ is zero. Thus, the rational points of $\C $ are the 
preimages of the rational torsion points of $\E$ and $\E'$. The torsion points 
of $\E$ are  $Tor (\E)=\{ \cO, (0,2), (0,-2) \}$. 
Then $\phi_1^{-1}(0)=\{0, \infty\}$ and 
 $\phi_1^{-1}(\infty)=\{-1, \frac 1 2 \pm \frac {\sqrt{-3}} 2 \}$.
Thus, $$\C(\bQ)=\{ (0,1), (0,-1), (-1,0)\}$$

\end{example}
%
%*********************************************************************
\begin{example} $\p=\p_6$:
 This point is in $\Y$ and it is not a singular point of $\Y$. 
It has isomorphic
degree 3 elliptic subcovers; see \cite{Sh2}. 
The corresponding $(\u, \v)$ pair is $(\u,\v)=(20, 16)$ and $e_3(\p)=1$. 
Then, the genus 2 curve has  equation:
$$\C: \quad Y^2=(256X^3+320X^2+16X+1 )\, (1024X^3+256X^2+32X+1)   $$
Degree 3 elliptic subcovers have  $j$-invariants $j=j'=-32768$ and equations 
\begin{equation}
\begin{split}
\E :&  \quad V_1^2= 4(- 5324 U_1^3 + 968 U_1^2 -56 U_1^2 +1 )\\
\E' :& \quad  V_2^2= 11^3 (-32000 \, U_2^3+35200 \, U_2^2-12320 \, U_2+11^3)\\
\end{split}
\end{equation}
where 
$$U_1= 16\, \frac {X^2} {256 X^3+320X^2+16X+1}, \quad U_2= \frac {1} {20} \,
\frac {(16X+3)^2 (20X+1)} {1024 X^3 + 256 X^2 + 32 X+1}$$
Both  elliptic curves have trivial torsion but rank $> 0$. One can try to adapt
more  sophisticated techniques in this case as Flynn and Wetherell have done for
degree 2 subcovers. 
This is the only genus 2 curve (up to $\bC$-isomorphism)
with  automorphism group  $D_{12}$ and isomorphic  degree 2 elliptic subcovers.
Indeed all degree 2 and 3 elliptic  subcovers are $\bC$-isomorphic with
 $ j$-invariants $-32768$. Degree 2 elliptic  subcovers also have rank 1 which does
not provide any quick information about rational points of $\C$. 
\end{example}
%*******************************************************************
%
\begin{example}
$\p=\p_7$:
All degree 2 and 3 elliptic subcovers are defined over $\bQ (\sqrt{5})$.
\end{example}

\medskip

\par Throughout this paper we have made use of several computer algebra packages as 
 {\sc Apecs},   {\sc Maple}, and {\sc GAP}. 
The interested reader can check  \cite{Sh-V} and 
\cite{Sh2} for more details on loci $\L_2$ and $\L_3$. 
The equations for these spaces, $j$-invariants of elliptic subcovers of degree 2 and 3,
and other computational aspects of genus 2 curves can be down loaded from author's 
web site.

%&&&&&&&&&&&&&&&&&&&&&&&&&&&&&&&&&&&&&&&&&&&&&&&&&&&&&&&&&&&&&&&&&&&


\begin{thebibliography}{99}

\bibitem {Bo} {\sc O. Bolza}, On binary sextics with linear transformations into
themselves. {\it Amer. J. Math.} {\bf 10}, 47-70.

\bibitem  {Cassels} {\sc J. Cassels  and V. Flynn}, Prolegomena to a Middlebrow
Arithmetic of Curves of Genus Two, LMS, 230, 1996.

\bibitem {Clebsch} {\sc A. Clebsch},
Theorie der Bin\"aren Algebraischen Formen, Verlag von B.G. Teubner, Leipzig, (1872).

\bibitem {ES}
{\sc T. Ekedahl, J. P.  Serre},
Exemples de courbes algébriques à jacobienne complètement décomposable. 
{\it C. R. Acad. Sci. Paris Sér. I Math.}, 317 (1993), no. 5, 509--513. 

\bibitem {FW} {\sc V. Flynn and J. Wetherell}, 
Finding rational points on bielliptic genus 2 curves, {\it Manuscripta Math.} 100,
519-533 (1999).

\bibitem {Fr} {\sc G. Frey},  On elliptic curves with isomorphic torsion
structures and corresponding curves of genus 2.  {\it Elliptic curves,
modular forms, and Fermat's last theorem (Hong Kong, 1993)}, 79-98,
Ser. Number Theory, I, {\it Internat. Press, Cambridge, MA}, (1995).

\bibitem {FK} {\sc G. Frey and E.  Kani},   Curves of genus 2 covering elliptic
curves and an arithmetic application. {\it Arithmetic algebraic
geometry (Texel, 1989)}, 153-176, {\it Progr. Math.}, 89, {\it
Birkh\"auser Boston, Boston, MA, (1991)}.

\bibitem {HLP} {\sc E. Howe, F. Lepr\'evost, and B. Poonen}, 
Large torsion subgroups of split Jacobians of curves of genus two or three. 
{\it Forum. Math}, {\bf 12} (2000), no. 3, 315-364.

\bibitem {Ig} {\sc J. Igusa}, Arithmetic Variety Moduli for genus 2. {\it
Ann. of Math}. (2), 72, 612-649, (1960).

\bibitem {KK}
{\sc W. Keller, L. Kulesz},
Courbes algébriques de genre 2 et 3 possédant de nombreux points rationnels.
 C. R. Acad. Sci. Paris Sér. I Math. 321 (1995), no. 11, 1469--1472. 

\bibitem {Krazer} {\sc A. Krazer}, Lehrbuch der Thetafunctionen, Chelsea, New York, (1970).

\bibitem {Ku} {\sc M. R. Kuhn}, Curves of genus 2 with split Jacobian. {\it Trans. Amer. Math. Soc}
{\bf 307}  (1988), 41-49

\bibitem {Lange} {\sc H. Lange}, {\"U}ber die Modulvariet\"at der Kurven vom
Geschlecht 2. {\it J. Reine Angew. Math.}, 281, 80-96, 1976.

\bibitem {MS} {\sc K. Magaard, T. Shaska, S. Shpectorov, and H. V\"olklein},
The automorphism group of a Riemann surface, 
 {\it  Proceedings of the Conference on
Communications in Arithmetic Fundamental Groups 
and Galois Theory} , Kyoto University,  Springer, 2001  (to appear).

\bibitem {Me} {\sc P. Mestre},
  Construction de courbes de genre 2 \'a partir de leurs modules. In T. Mora 
and C. Traverso, editors,
 {\it Effective methods in algebraic geometry}, volume 94. {\it Prog. Math.
}, 313-334. Birkh\"auser, 1991. Proc. Congress in Livorno, Italy, April 17-21, (1990).

\bibitem {Mu} {\sc D. Mumford}, The Red Book of Varieties and Schemes, Springer, 1999.

\bibitem {Sh1}{\sc T.  Shaska}, Genus 2 curves with (n,n)-decomposable
 Jacobians, {\it Jour. Symb. Comp.}, Vol 31, {\bf no. 5}, pg.
603-617, 2001.

\bibitem {Sh2} {\sc T. Shaska}, Genus 2 fields with degree 3 elliptic 
subfields, (submited for publication).

\bibitem {Sh-V} {\sc T. Shaska and H.  V\"olklein}, 
Elliptic Subfields and automorphisms
of genus 2 function fields.  {\it Proceeding of the Conference on
Algebra and Algebraic Geometry with Applications: The celebration
of the seventieth birthday of Professor S.S. Abhyankar},
Springer-Verlag, 2001.

\end{thebibliography}
\end{document}